\documentclass[11pt]{article}
\usepackage[T2A]{fontenc}
\usepackage[cp1251]{inputenc}
\usepackage[russian]{babel}
\usepackage{amsmath,amsfonts,amssymb,amsthm}
\usepackage{ifthen}

\advance\textwidth40mm
\advance\textheight52mm
\advance\voffset-33mm
\advance\hoffset-25mm

\newtheorem{thm}{Теорема}[section]
\newtheorem{cor}[thm]{Следствие}
\newtheorem{df}[thm]{Определение}
\newtheorem{lm}[thm]{Лемма}
\newtheorem{prop}[thm]{Предложение}

\makeatletter

\newcommand\txt[2]{{#1{\ifthenelse{\equal {#2}{}}{}{ #2}}. }}
\def\note{\par {\bf Замечание.} }
\def\example#1.{\txt{\bf Пример}{#1}}
\makeatother

\def\itref#1{{\it (#1)}}

\newcommand\emp\varnothing
\newcommand\eps\varepsilon
\newcommand\ov\widetilde
\let\ds\displaystyle

\newcommand\eop{\unskip\hfill\penalty150\hfill$\square$\par\smallskip}

\def\SC{{\cal S}}
\def\OC{{\cal O}}
\def\II{{\cal I}}

\begin{document}

\title{Наибольшая длина периода слова, задаваемого $n$ запретами}
\author{И.И.~Богданов, Г.Р.~Челноков}

\maketitle

\footnotetext[1]{The work is supported by the Russian government project 11.G34.31.0053.}

\section{Введение}

Исследование комбинаторных свойств периодических последовательностей
(слов) играет важную роль в проблемах бернсайдовского типа,
см., например, \cite{Ufn,Kur,BBL}.

Алгебра $A=F\langle x_1,\dots,x_n\rangle/I$ называется {\em мономиальной}, если идеал~$I$ свободной алгебры~$F\langle x_1,\dots,x_n\rangle$ порождён мономами. При изучении мономиальных алгебр важную роль играют алгебры вида~$A_u$, где $u$ --- непериодичное слово: алгебра~$A_u$ задана соотношениями $v=0$, где $v$ пробегает множество всех слов, не являющихся подсловами в~$u^\infty$, см.~\cite{BBL}. Алгебрами~$A_u$ исчерпывается класс первичных конечно определенных
мономиальных PI-алгебр.
В то же время, не все слова~$v$
необходимы для задания такой алгебры. Достаточно, например,
ограничиться всеми словами длины, не превосходящей~длины $u$.

Представляет интерес более точное исследование структуры соотношений, задающих
$A_u$. В данной работе исследуется вопрос о возможной длине слова~$u$, при
котором алгебра~$A_u$ может быть задана $n$ мономиальными соотношениями. Мы показываем (см.
теорему~\ref{main}), что в случае алфавита из двух букв наибольшая длина слова равен
числу Фибоначчи~$F(n)$.

Работа является продолжением статьи~\cite{cheln}, в которой получены экспоненциальные оценки на длину слова~$u$. Мы используем некоторые понятия и результаты из этой статьи.

Как авторам стало известно, в настоящее время П. Лавров предложил другое доказательство этого факта~\cite{lavrov}. Было бы интересно сравнить методы доказательств.

\section{Предварительные сведения}

Пусть~$X=\{x_1,\dots,x_k\}$ --- конечный алфавит (в большей части статьи мы
полагаем $k=2$). Под {\em конечным (бесконечным вправо/влево/в обе стороны) словом} мы понимаем любую конечную (бесконечную вправо/влево/в обе стороны) последовательность букв алфавита; пустая
последовательность~$\Lambda$ также является словом. {\em Длиной}~$|u|$
конечного слова~$u$ называется количество букв в нём.
Все конечные слова образуют моноид относительно конкатенации.

\begin{df}
  Слово $u$ называется {\em подсловом} слова~$w$, если $w=v_1uv_2$ для некоторых слов~$v_1$, $v_2$. Слово $u$ является {\em началом (концом)} слова~$w$, если $v_1=\Lambda$ ($v_2=\Lambda$). Подслово (начало, конец)~$u$ слова~$v$ является {\em собственным}, если $u\neq v$.

  Введём на множестве конечных слов частичный порядок: скажем, что $u\preceq v$, если $u$ является подсловом слова~$v$.

  Непустое слово $u$ называется {\em периодическим}, если $u=v^n$ для некоторого слова $v$ и некоторого $n\geq 2$. В противном случае оно называется {\em непериодическим}.
\end{df}

\begin{df}[\cite{cheln}]
  {\em Системой запретов} назовём конечное множество слов $V=\{v_1,\dots,v_n\}$ в алфавите~$X$.
  Будем говорить, что (конечное или бесконечное) слово~$w$ {\em удовлетворяет системе запретов} $V$, если $v\not\preceq w$ для любого $v\in V$.

  Пусть~$W$ --- бесконечное в обе стороны слово. Будем говорит, что система
  запретов~$V$ {\em определяет} слово~$W$, если $W$ --- единственное бесконечное слово,
  удовлетворяющее этой системе запретов.
\end{df}

Пусть $u$ --- конечное слово. Определим {\em бесконечное слово с периодом $u$}
как $u^\infty=\dots uuuu\dots$. Если существует такое слово~$u$, что $W=u^\infty$, то бесконечное слово~$W$ назовём {\em периодичным}. Нетрудно видеть, что если система запретов определяет слово~$W$, то оно периодично.

Для каждого периодичного бесконечного в обе стороны слова~$W$ существует в
определённом смысле оптимальная система запретов. Слово $v$ назовём {\em
каноническим запретом} для~$W$, если $v$ не является подсловом~$W$, а любое его
собственное подслово --- является. Множество всех канонических запретов
для~$W$ назовём {\em канонической системой запретов} для~$W$; она
обозначается~$C(W)$.

\begin{lm}[{см. \cite[Лемма 1]{cheln}}]
  \label{min-szap}
  Каноническая система запретов~$C(W)$ определяет слово~$W$. При этом
  любая система запретов, задающая~$W$, содержит не меньше элементов,
  чем~$C(W)$.
  \qed
\end{lm}

\note Можно показать также, что $C(W)$ --- единственная система запретов,
задающая~$W$, с минимальной возможной суммой длин входящих в неё слов.

Отметим ещё одно полезное свойство системы~$C(W)$.

\begin{prop}
  Конечное слово $v$ удовлетворяет~$C(W)$ тогда и только тогда, когда $v$ --- подслово
  в~$W$.
\end{prop}

\proof
Если $v$ --- подслово в~$W$, то, очевидно, оно удовлетворяет~$C(W)$. Обратно,
предположим, что $v$ не является подсловом в~$W$. Тогда существует минимальное
по длине подслово $v'\preceq v$, не являющееся подсловом~$W$; оно по
определению лежит в~$C(W)$. Значит, $v$ не удовлетворяет~$C(W)$.
\eop

Определим числа Фибоначчи по следующему правилу: $F(0)=F(1)=1$, $F(k+1)=F(k)+F(k-1)$. Мы продолжим эту последовательность на отрицательные индексы; так, $F(-1)=0$, $F(-2)=1$, $F(-3)=-1$.

Цель данной статьи --- нахождение точных верхних оценок на
длину периода слова, если известна мощность системы запретов, его
задающая. Основным результатом статьи является следующая теорема.

\begin{thm}
  \label{main}
  Пусть $|X|=2$. Пусть система запретов~$V$ определяет слово~$W=u^\infty$, где слово $u$ непериодично. Тогда $|u|\leq F(|V|)$.
\end{thm}

\note В силу леммы~\ref{min-szap} можно ограничиться случаем $V=C(W)$.

\begin{df}
  Пусть слово~$u$ удовлетворяет системе запретов~$V$, и $x\in X$. Назовём
  слово~$u'=ux$ ($u'=xu$) {\em продолжением} слова~$u$ {\em вправо (влево)} относительно~$V$,
  если $u'$ также удовлетворяет~$V$.

  Слово~$u$ назовём {\em неоднозначно продолжимым вправо
  (влево)}, если у него существуют хотя бы два разных продолжения вправо
  (влево).

  Наконец, назовём~$u$ {\em развилкой} (относительно~$V$), если $u$
  неоднозначно продолжимо как вправо, так и влево. {\em Кратностью}
  развилки~$u$ назовём количество её продолжений вправо.

  Назовём слово~$u$ {\em развилкой} относительно бесконечного в обе стороны
  слова~$W$, если~$u$ является развилкой относительно~$C(W)$. Само слово~$W$
  также назовём {\em развилкой} относительно~$W$.
\end{df}

\example. Пусть $|X|=2$, $X=\{a,b\}$. Тогда конечное слово~$u$ является развилкой
относительно~$W$ тогда и только тогда, когда все четыре слова $ua$, $ub$, $au$,
$bu$ являются подсловами в~$W$. При этом все развилки имеют кратность~2. Стоит
отметить, что слова $aua$, $aub$, $bua$, $bub$ уже не обязательно являются
подсловами в~$W$; с другой стороны, нетрудно видеть, что хотя бы два из них
должны удовлетворять~$C(W)$.

\smallskip
В работе~\cite{cheln} задача оценки количества запретов, задающих слово~$W$, была сведена к задаче оценки количества развилок в слове~$W$. Мы также будем использовать этот результат.

\begin{lm}[{см.~\cite[Лемма 3]{cheln}}]
  Для каждого подслова~$u$ слова~$W$ существует наименьшая (относительно порядка $\preceq$) развилка~$v=r(u)$, содержащая~$u$.
  \qed
\end{lm}

\note Если $v\preceq v'$ --- подслова в~$W$, то, очевидно, $r(v)\preceq r(v')$.

\medskip

\begin{lm}[{см.~\cite[Лемма 5]{cheln}}]\
  \label{kratn-zp}
  Пусть $v_1,\dots,v_n$ --- все конечные развилки в периодичном слове~$W=u^\infty$, а $k_1,\dots,k_n$ --- их кратности. Тогда
  $$
    |C(W)|\geq 1+(k_1-1)+(k_2-1)+\dots+(k_n-1).
  $$
\end{lm}

\begin{cor}
  \label{cor-kratn-zp}
  $|C(W)|\geq n+1$.
  \eop
\end{cor}

\note Можно показать, что при $|X|=2$ в лемме~\ref{kratn-zp} и в следствии~\ref{cor-kratn-zp}
всегда достигается равенство.

\smallskip
Это следствие позволяет свести теорему~\ref{main} к следующей.

\begin{thm}
  \label{main'}
  Пусть $|X|=2$. Рассмотрим периодическое слово~$W=u^\infty$, где слово $u$ непериодично. Пусть для этого слова существует $n$
  конечных развилок. Тогда $|u|\leq F(n+1)$.
\end{thm}

Именно этот вариант мы и доказываем в конце раздела~\ref{bound}.

\medskip

В заключение приведём пример, показывающий, что оценка в теореме~\ref{main} (и,
следовательно, в теореме~\ref{main'}) неулучшаемы ни при каком $n\geq 2$.

\example. Пусть $X=\{a,b\}$. Построим последовательности слов $(s_i)$, $(t_i)$
по следующему правилу. Положим $s_0=a$, $t_0=b$; далее, при всех $i\geq 0$
положим $s_{i+1}=s_is_it_i$, $t_{i+1}=s_it_i$. Нетрудно видеть, что
$|s_i|=F(2i+1)$, $|t_i|=F(2i)$. В работе~\cite[Теорема~3]{cheln} показано, что
при $i\geq 1$ слово $W_{2i}=(t_i)^\infty$ задаётся $2i$ запретами; значит,
$|C(W_{2i})|\leq 2i$ по лемме~\ref{min-szap}. Тогда ясно, что слово $W_{2i}$
показывает неулучшаемость оценок в теоремах~\ref{main} и~\ref{main'} при
чётном~$n$.

Аналогично можно показать, что при $i\geq 1$ слово $W_{2i+1}=(s_i)^\infty$
задаётся $2i+1$ запретом; это показывает неулучшаемость оценок при
нечётном~$n$.

\section{Комбинаторика}

На протяжении этого и последующего разделов мы рассматриваем фиксированное непустое конечное непериодическое слово~$u$, и слово $W=u^\infty$. Развилки и запреты относительно слова~$W$ мы называем просто развилками и запретами.

\begin{df}
  Назовем {\em значимостью} $z(v)$ подслова~$v$ количество раз, которое оно встречается на периоде;
  формально говоря, если $u=u_1\dots u_d$, где $u_1,\dots,u_d\in X$, и $|v|=t$, то
  $$
    z(v)=\left|\left\{1\leq i\leq d: u_i\dots u_{i+t-1}=v\right\}\right|,
  $$
  где мы полагаем $u_{t+d}=u_i$ при $1\leq i\leq d$.
\end{df}

Напомним, что для подслова~$v$ слова~$W$ через $r(v)$ обозначается наименьшая
развилка, содержащая~$v$.

\begin{prop}
  Если $v\preceq v'$, то $z(v)\leq z(v')$. Кроме того, $z(v)=z(r(v))$.
  \label{leqzn}
  \qed
\end{prop}

%
\begin{prop}
  Пусть $v$~--- произвольная конечная развилка. Тогда
  $$
    z(v)=\sum_{x\in X}z(vx)=\sum_{x\in X}z(r(vx)).
    \eqno{\square}
  $$
  \label{sumzn}
\end{prop}


Пусть $v_0,v_1,\dots,v_n$~--- все развилки, упорядоченные по значимости, то
есть $z_0\leq z_1\leq\dots\leq z_n$, где $z_i=z(v_i)$. При этом мы считаем, что
$v_0=W$ (и $z_0=1$), а $v_n=\Lambda$ (и $z_n=|u|$). Таким образом, наша цель~--- получить верхнюю оценку на~$z_n$.

Из предложения~\ref{sumzn} следует, что $z_1=z_0+z_0=2$. Из предложений~\ref{leqzn} и~\ref{sumzn} следует следующее предложение.

\begin{prop}
  Пусть $x\in X$, $0\leq i\leq n$. Тогда $z(v_ix)<z(v_i)$. В частности, если $r(v_ix)=v_j$, то
  $j<i$. Наконец, из $v_i\prec v_j$ следует $z_i>z_j$.
  \eop
\end{prop}

Далее мы работаем со словами в алфавите $X=\{a,b\}$. В этом случае кратность
каждой развилки равна~2. Из предложения~\ref{sumzn} теперь вытекает следующее предложение.

\begin{prop}
  $z_i\leq 2z_{i-1}$, и $\max\{r(v_ia), r(v_ib)\} \ge z_i/2$.
  \eop
  \label{dbling}
\end{prop}

\begin{df}
  Назовем развилку $v_i$ ($i\geq 2$) {\em исключительной},
  если $z_i>z_{i-1}+z_{i-2}$. В противном случае назовем $v_i$ {\em регулярной}.
  Развилки $v_0$ и $v_1$ также будем считать регулярными. Индекс $i$ назовем {\em
  исключительным (регулярным)}, если развилка~$v_i$ исключительна (регулярна). Обозначим множество исключительных развилок через~$\II$.
\end{df}

Неформально говоря, в регулярных случаях последовательность $(z_i)$ растет не
быстрее чисел Фибоначчи.

\begin{prop}
  Если развилка $v_i$ исключительна, то $z_i=2z_{i-1}$, $z_{i-1}>z_{i-2}$, и $r(v_ia)=r(v_ib)=v_{i-1}$.
  \label{regcase}
\end{prop}

\proof
  Пусть $r(v_ia)=r(v_ib)=v_{i-1}$; тогда $z_i=2z_{i-1}$, и исключительность
  развилки~$v_i$ равносильна тому, что $2z_{i-1}>z_{i-1}+z_{i-2}$, то есть
  $z_{i-1}>z_{i-2}$. В противном случае можно считать, что $r(v_ia)=v_j$ при
  $j\leq i-2$. Тогда по предложению~\ref{sumzn}
  $z_i=z(v_i)=z(r(v_ia))+z(r(v_ib))\leq z_{i-2}+z_{i-1}$, то есть $v_i$ регулярна.
\eop

\note Исключительные развилки могут существовать. Например, в слове $u^\infty$, где $u=(ababbabbabbb)^na$, развилка $v=babbabb$ исключительна при $n\geq 2$. Действительно, нетрудно проверить, что $z(v)=2n$, $r(va)=r(vb)=ababbabbabbba=w$, $z(w)=n$; значимость же любой другой развилки либо не меньше $3n-1$, либо не больше $n-1$.

\smallskip

Остаток этого раздела посвящён изучению исключительных развилок.

\begin{df}
  \label{straff}
  Пусть $v_i\in\II$. Пусть $v_j$ --- максимальное собственное начало развилки $v_{i-1}$, являющееся развилкой. Назовём развилку $v_j$ и её индекс~$j$ {\em штрафными} для исключительной развилки~$v_i$ и её индекса~$i$; мы будем обозначать $v_j=\Psi(v_i)$.
\end{df}

\note В принципе, определением не запрещена ситуация $i=j$; но в дальнейшем мы увидим, что она невозможна, см. предложение~\ref{est-pen}.


Из предложения~\ref{sumzn} вытекает

\begin{prop}
  Пусть $v_i\in\II$ и $v_j=\Psi(v_i)$. Тогда $z_j\leq z_{j-1}+z_{i-1}$.
  \label{est-straf}
  \eop
\end{prop}

\begin{prop}
    Пусть $v_i\in\II$ и $v_j=\Psi(v_i)$, причём $v_{i-1}=r(v_ja)$. Тогда $v_{i-1}\succeq v_jb$.
  \label{sb in q}
\end{prop}

\proof Построим последовательность развилок $(s_k)$ следующим
образом. Положим $s_0=v_j$, $s_1=r(v_jb)$; заметим, что
$z(s_1)=z(r(v_jb))=z(v_j)-z(r(v_ja))=z(v_j)-z(v_{i-1})\geq
z(v_i)/2$, так как $z(v_j)\geq z(v_i)=2z(v_{i-1})$. При $k\geq 1$
через $s_{k+1}$ обозначим такую из развилок~$r(s_ka)$ и~$r(s_kb)$,
для которой~$z(s_{k+1})\geq z(s_k)/2$; она существует согласно
предложению~\ref{dbling} (по замечанию выше, неравенство
$z(s_{k+1})\geq z(s_k)/2$ выполнено и при $k=0$). Заметим, что
$v_jb\preceq s_k$ при каждом~$k\geq 1$.

Пусть $k$~--- максимальное число, для которого $z(s_k)\geq
z(v_{i-1})$; пусть $s_k=v_m$. Предположим, что $m\ne i-1$. По
предложению~\ref{regcase} имеем $z_{i-2}<z_{i-1}$; значит, случай
$m<i-1$ невозможен. Поэтому $m\geq i$, то есть $z(s_k)\geq
z_i=2z_{i-1}$. Но тогда $z(s_{k+1})\geq z(s_k)/2\geq z_{i-1}$, что
противоречит выбору~$k$. Итак, $m=i-1$, поэтому $v_{i-1}=s_k\succeq
v_jb$. \eop

\begin{prop}
  Пусть $v_i\in\II$ и $v_j=\Psi(v_i)$, причём $v_{i-1}=r(v_ja)$. Тогда существует такое $k$ ($i<k<j$), что $z_k\leq z_{k-1}+z_{i-2}$ и $z_k<z_j$. В частности, $j\geq i+2$, и развилка $v_k$ регулярна. Кроме того, в канонической системе запретов существует запрет вида $yv_ka$, где $y$ "--- буква.
  \label{est-pen}
\end{prop}

\proof
Согласно определению слова~$v_j$ и предложению~\ref{sb in q}, слово $v_{i-1}$ можно представить как $v_{i-1}=v_jat_1=t_2v_jbt_3$
для некоторых слов $t_1$, $t_2$, $t_3$. Слово $t_2$, очевидно,
непусто; пусть $x$~--- его последняя буква, $t_2=t_2'x$. Поскольку
$v_{i-1}=r(v_ja)$, любое вхождение $v_ja$ в слово~$W$ продолжается
до~$v_jat_1=v_{i-1}$; в частности, оно продолжается до $t_2v_jb$.
Это значит, что слово $t_2v_ja$ (начинающееся с~$v_ja$) не
встречается в~$W$.

Тогда слово $t_2v_ja=t_2'xv_ja$ должно содержать некоторый запрет
$yv_kz$ из канонической системы (здесь $y$, $z$~--- буквы, $v_k$~---
некоторая развилка). Этот запрет не может быть подсловом
слова~$t_2v_j$, ибо оно встречается в~$W$. Также он не может являться
подсловом слова~$xv_ja$. Действительно, поскольку $v_{i-1}$ является
развилкой, слова $av_{i-1}$ и $bv_{i-1}$ встречаются в~$W$; значит,
и их подслова $av_ja$ и $bv_ja$ также в нем встречаются и потому не
могут содержать запретов.

Итак, наш запрет $yv_kz$ не содержится в подсловах $t_2v_j$ и
$xv_ja$. Это значит, что он является концом слова $t_2v_ja$, строго
содержащим~$xv_ja$; таким образом, $z=a$, а $v_k=s'v_j$ для
некоторого непустого слова~$s'$. Рассмотрим теперь развилку~$v_\ell=r(v_ka)$. Слово $v_ka$ заканчивается на $v_ja$; значит, развилка $v_\ell$ должна содержать развилку $r(v_ja)=v_{i-1}$. Более того, согласно определению, слово $v_ja$ является началом развилки~$v_{i-1}$ и находится не в начале развилки~$v_\ell$; значит, $v_{i-1}$ "--- собственное подслово в~$v_\ell$, то есть $v_\ell\succ
v_{i-1}$. Поскольку и $v_\ell$, и $v_{i-1}$ "--- развилки, получаем, что $z(v_\ell)<z(v_{i-1})$ и $\ell\leq i-2$.

Итого, мы нашли развилку $v_k=s'v_j$ такую, что $z(r(v_ka))\leq
z_{i-2}$; значит, $z_k=z(v_k)=z(r(v_ka))+z(r(v_kb))\leq
z_{i-2}+z_{k-1}$. Заметим, что $v_j\prec v_k\prec v_{i-1}$,
поэтому $i-1<k<j$ и $z_k=z(v_k)<z(v_j)=z_j$. Кроме того, $k\ne i$,
ибо $z(r(v_ia))=z_{i-1}>z_{i-2}\geq z_\ell=z(r(v_ka))$. Значит,
$i<k<j$ (и, значит, $j\geq i+2$), и требуемое~$k$ найдено. Наконец,
поскольку $z_k\leq z_{k-1}+z_{i-2}<z_{k-1}+z_{i-1}\leq
z_{k-1}+z_{k-2}$, развилка $v_k$ регулярна.
\eop

\begin{df}
  Пусть $v_i\in \II$, $v_j=\Psi(v_i)$. Пусть $v_k$ --- развилка, построенная в предложении~\ref{est-pen}. Назовем эту развилку~$v_k$ и ее индекс~$k$ {\em пеневыми} для исключительной развилки~$v_i$ и ее индекса~$i$; обозначим $v_k=\Pi(v_i)$.
\end{df}

Отметим некоторые свойства штрафных и пеневых развилок.

\begin{prop}
  Пусть $v_i\in\II$, $v_j=\Psi(v_i)$ и $v_{i-1}=r(v_ja)$. Тогда $z(r(v_ja))<z(r(v_jb))$. В частности, развилка~$v_j$ регулярна.
  \label{str-reg}
\end{prop}

\proof Первое утверждение следует из того, что $r(v_ja)=v_{i-1}$, а
$z(r(v_ja))+z(r(v_jb))=z(v_j)=z_j>z_k\geq z_i=2z(r(v_ja))$. Тогда
$v_j$ не исключительна по предложению~\ref{regcase}. \eop

\begin{prop}
  Пусть $v_i,v_{i'}\in \II$, $v_j=\Psi(v_i)$, $v_{j'}=\Psi(v_{i'})$. Тогда, если $i\neq i'$, то и $j\neq j'$.
  \label{dist-straf}
\end{prop}

\proof
Предположим противное; пусть $i>i'$ и $v_{i-1}=r(v_ja)$. Тогда из
предложений~\ref{regcase} и~\ref{str-reg} следует, что $z_{i'-1}<
z_{i-1}=z(r(v_ja))<z(r(v_jb))$, и потому $v_{i'-1}$ не может являться
$r(v_jb)$. Таким образом, $v_{i'-1}=r(v_ja)=v_{i-1}$, и $i=i'$. Противоречие.
\eop

\begin{prop}
  Пусть $v_i,v_{i'}\in \II$, $v_k=\Pi(v_i)$. Тогда $v_k\neq \Psi(v_{i'})$.
  \label{dist-pen}
\end{prop}

\proof
Пусть $v_j=\Psi(v_i)$, причём $v_{i-1}=r(v_ja)$. По предложению~\ref{est-pen}),
в канонической системе запретов существует запрет вида~$yv_ka$, где $y$~--- буква, при этом $z(r(v_ka))\leq z_{i-2}$, а $z(v_k)\geq z_i> 2z_{i-2}$. Значит, $z(r(v_kb))=z(v_k)-z(r(v_ka))> z(r(v_ka))$. Поэтому, если развилка~$v_k=\Psi(v_{i'})$, то
по предложению~\ref{str-reg} $v_{i'-1}=r(v_ka)$, и $v_ka$ является началом
слова~$v_{i'-1}$ по определению штрафной развилки. Но поскольку~$v_{i'-1}$
является развилкой, то подслово $yv_{i'-1}$ (и тем более $yv_ka$) встречается
в~$W$ и потому не может являться запретом~--- противоречие.
\eop

Суммируем результаты предложений~\ref{est-straf}, \ref{est-pen}, \ref{str-reg},
\ref{dist-straf} и~\ref{dist-pen} в следующей теореме.

\begin{thm}
  \label{struct}
  Для каждого исключительного индекса~$i$ существуют штрафной и пеневой
  индексы~$j$ и~$k$ такие, что $i<k<j$,
  $z_j\leq z_{j-1}+z_{i-1}$ и $z_{k}\leq z_{k-1}+z_{i-2}$.
  При этом штрафные индексы для разных исключительных также различны,
  а пеневой не может являться штрафным. Кроме того, штрафные и пеневые индексы
  регулярны (т.е. не исключительны).
  \eop
\end{thm}

\begin{df}
  Назовем индекс $r$ {\em рядовым}, если он не
  является ни исключительным, ни штрафным, ни пеневым.
\end{df}

\section{Оценки}
\label{bound}

В этом разделе мы оцениваем рост последовательности~$(z_i)$. Для этого мы сначала введём класс абстрактных (не обязательно связанных со словами) последовательностей, мажорирующих последовательности вида $(z_i)$, а затем будем оценивать эти последовательности.

\subsection{Допустимые последовательности}

\begin{df}
  Пусть $n\geq 2$ --- натуральное число. Пусть в множестве $\{2,3,\dots,n\}$ выделены три попарно
  непересекающихся подмножества $I$, $J$ и $K$, $|I|=|J|\geq |K|$ (элементы этих подмножеств
  будем называть соответственно {\em исключительными, штрафными и пеневыми};
  индекс, не лежащий ни в одном из подмножеств, назовем {\em рядовым}).
  Наконец, пусть вдобавок зафиксированы биекция $\psi:I\to J$ и
  сюръекция $\pi: I\to K$, причем $i<\pi(i)<\psi(i)$
  для любого $i\in I$. Назовём набор $\SC=(n,I,J,K,\psi,\pi)$ {\em системой}.
  Для $k\in K$ определим $d(k)=\min \pi^{-1}(k)$; элементы множества
  $d(K)\subseteq I$ назовём {\em плохими} для системы~$\SC$.

  Самой простой системой является <<пустая>> система $\OC_n=(n,\emp,\emp,\emp,\emp,\emp)$.

  Пусть $\Pi=(x_i)_{i=0}^n$ --- последовательность неотрицательных чисел. Будем говорить, что~$\Pi$ {\em соответствует} системе~$\SC$, если выполнено условие

  (1)  для любого $2\leq r\leq n$,
  $x_r=x_{r-1}+x_{\theta(r)}$, где
  $$
    \theta(r)=
    \begin{cases}
      r-2, & \text{если $r$ --- рядовой;}\\
      r-1, & \text{если $r$ --- исключительный;}\\
      \psi^{-1}(r)-1, & \text{если $r$ --- штрафной;}\\
      d(r)-2, & \text{если $r$ --- пеневой.}
    \end{cases}
  $$
  Ясно, что такая последовательность задаётся начальными членами $x_0$ и $x_1$; будем обозначать её $\Pi_\SC(x_0,x_1)$.

  Назовём последовательность~$\Pi_\SC(a,b)$ {\em допустимой} для~$\SC$, если $0\leq a\leq b\leq 2a$; наконец, будем говорить, что последовательность $\Pi_\SC=\Pi_\SC(1,2)$ {\em порождена} системой~$\SC$.
\end{df}

Пусть теперь $W=u^\infty$ --- бесконечное периодичное слово, и $(z_i)$ --- последовательность значимостей, определённая в предыдущем разделе. Результаты этого раздела позволяют выписать порождённую последовательность, мажорирующую последовательность~$(z_i)$. Именно,
пусть индексы $i_1<\dots<i_m$ являются исключительными для слова~$W$, индексы $j_1,\dots,j_m$~--- штрафными (причем $j_t$~--- штрафной для $i_t$), а индексы $k_1,\dots,k_s$~---
пеневыми (напомним, что пеневой индекс может соответствовать
нескольким исключительным). Положим $I=\{i_1,\dots,i_m\}$,
$J=\{j_1,\dots,j_m\}$, $K=\{k_1,\dots,k_s\}$; по
теореме~\ref{struct} эти множества попарно не пересекаются. Далее,
для всех $1\leq r\leq m$ положим $\psi(i_r)=j_r$ и определим
$\pi(i_r)$ как пеневой индекс, соответствующий
исключительному~$i_r$. Тогда $(n,I,J,K,\psi,\pi)$
--- система согласно теореме~\ref{struct}. Из этой же теоремы
вытекает следующее предложение.

\begin{prop}
  \label{z<y}
  Пусть последовательность $(y_i)$ порождена системой
  $\SC=(n,I,J,K,\psi,\pi)$. Тогда для любого индекса~$r=0,\dots,n$
  выполнено неравенство $z_r\leq y_r$.
\end{prop}

\proof
Индукция по~$r$. При $r=0,1$ утверждение очевидно, так как $z_r=y_r$. Пусть
$z_s\leq y_s$ при всех $s<r$. Тогда по теореме~\ref{struct} имеем $z_r\leq
z_{r-1}+z_{\theta(r)}\leq y_{r-1}+y_{\theta(r)}=y_r$, что и требовалось.
\eop

Отметим сразу некоторые свойства любой допустимой последовательности $(y_i)$, аналогичные свойствам последовательности~$(z_i)$ из предыдущего раздела.

\begin{prop}
  \label{vdvoe}
  Пусть последовательность~$(y_i)$ соответствует системе~$\SC$. Тогда $0\leq y_i\leq y_{i+1}$ при всех $1\leq i\leq n$. Если, вдобавок, $(y_i)$ допустима для~$\SC$, то $0\leq y_i\leq y_{i+1}\leq 2y_i$ при всех $0\leq i<n$.
\end{prop}

\proof
Неравенство $y_i\geq 0$ (и поэтому $y_{i+1}\geq y_i$) следует из определения. Осталось доказать неравенство~$y_i\leq y_{i+1}\leq 2y_i$ для допустимой последовательности~$(y_i)$. Применим индукцию по~$i$. При $i=0$ все утверждения верны. Далее, при $i\geq 1$
имеем $y_{i+1}=y_i+y_r$ при некотором $r=\theta(i+1)\leq i$. По предположению
индукции имеем $0\leq y_r\leq y_i$, откуда $y_i\leq y_i+y_r\leq 2y_i$, что и
требовалось доказать.
\eop

Напомним, что числа Фибоначчи заданы условиями $F(0)=F(1)=1$ и $F(n+1)=F(n)+F(n-1)$ при всех целых~$n$.

\begin{prop}
  \label{ryad}
  Пусть $k\geq 2$, $t\geq -1$, и индексы $k,k+1,\dots,k+t$ --- рядовые. Тогда
  $y_{k+t}=F(t+1)y_{k-1}+F(t)y_{k-2}$.
\end{prop}

\proof
Индукция по $t$. При $t=-1,0,1$ имеем
$$
  y_{k-1}=F(0)y_{k-1}+F(-1)y_{k-2}, \quad y_k=F(1)y_{k-1}+F(0)y_{k-2}, \quad y_{k+1}=y_{k-1}+y_k=F(2)y_{k-1}+F(1)y_{k-2}.
$$
Если же $t\geq 2$, то
$$
  y_{k+t}=y_{k+t-1}+y_{k+t-2}=(F(t)+F(t-1))y_{k-1}+(F(t-1)+F(t-2))y_{k-2}
    =F(t+1)y_{k-1}+F(t)y_{k-2},
$$
что и требовалось.
\eop

\medskip

Пусть последовательность $(x_i)$ порождена системой~$\OC_n$; тогда, ясно, $x_i=F(i+1)$ при всех $0\leq i\leq n$. Наша цель --- показать, что для любой порождённой последовательности $y_0,\dots,y_n$ выполняется неравенство $y_n\leq x_n=F(n+1)$. С этой целью мы будем перестраивать
систему $(n,I,J,K,\psi,\pi)$, сводя её к пустой, так, чтобы значение $y_n$ не
уменьшалось.

\subsection{Элементарные улучшения порождённой последовательности}

Здесь и далее, если не оговорено противное, $\SC=(n,I,J,K,\psi,\pi)$ ---
произвольная система, а $y_0,\dots,y_n$ "--- последовательность, ею порождённая. Обозначим через $L=I\cup
J\cup K$ множество всех нерядовых индексов.

Каждый раз мы будем перестраивать систему $\SC$, получая систему $\SC'=(I',J',K',\psi',\pi')$ и порождённую ею последовательность $(y_i')_{i=0}^n$, для которой $y_n'\geq y_n$ (функции $d$ и $\theta$, а также множество~$L$ для системы~$\SC'$ также будем помечать штрихами). Такую последовательность $(y_i')$ (систему $\SC'$) мы будем называть {\em улучшением} последовательности $(y_i)$ (системы~$\SC$). Достаточные условия для улучшения обеспечивает следующая лемма.

\begin{lm}[об улучшении]
  \label{uluchsh}
  Пусть $\ell\geq 2$, и выполнены следующие условия:

%
  (1) из $\theta(i)\geq\ell-1$ следует $\theta'(i)=\theta(i)$;

  (2) $y_{\ell-1}'\geq y_{\ell-1}$, $y_{\ell}'\geq y_{\ell}$;

  (3) $y_{\theta'(i)}'\geq y_{\theta(i)}$ при всех $i>\ell$ таких, что $\theta(i)<\ell-1$.

  Тогда $y_i'\geq y_i$ при всех $i\geq \ell-1$; в частности, $(y_i')$ является
  улучшением~$(y_i)$.
\end{lm}

\proof
Индукция по~$i$. База для $i=\ell-1,\ell$ есть условие~\itref{2}; пусть $i>\ell$. Если $\theta(i)\geq \ell-1$, то $y_i=y_{i-1}+y_{\theta(i)}\leq y_{i-1}'+y_{\theta'(i)}'=y_i'$ по предположению индукции и условию~\itref{1}. Если же $\theta(i)<\ell-1$, то $y_i=y_{i-1}+y_{\theta(i)}\leq y'_{i-1}+y'_{\theta'(i)}=y'(i)$ по предположению индукции и условию~\itref{3}.
\eop

\smallskip

В этом подразделе мы приведём несколько элементарных улучшений. Первые два из
них можно схематично изобразить так:
$$
  \text{ИР}\dots\text{РН}\to\text{Р}\dots\text{РИН};
  \qquad\qquad
  \text{ИКР}\to\text{КРИ}, \quad \text{ИКН}\to \text{КИН},
$$
где через И, Р, К, Н обозначены соответственно исключительный индекс, рядовой
индекс, штрафной или пеневой индекс, нерядовой неплохой индекс (Напомним, что индекс плох, если он лежит в множестве~$d^{-1}(K)$.

\begin{prop}[сдвиг исключительного индекса вправо]
  \label{iskl-right}
  Пусть $r$ --- исключительный индекс, а $\ell=\min\{t:r<t\in L\}$ --- следующий за $r$
  нерядовой индекс. Пусть индекс $\ell$ неплохой. Обозначим через $I'$ множество, полученное из $I$ заменой $r$ на
  $r'=\ell-1$. Соответственно изменим функции $\psi$, $\pi$, полагая
  $\psi'(r')=\psi(r)$, $\pi'(r')=\pi(r)$. Тогда система
  $\SC'=(n,I',J,K,\psi',\pi')$ --- улучшение системы~$\SC$.
\end{prop}

\proof
При $\ell=r+1$ доказывать нечего; пусть $\ell\geq r+2$. Очевидно, что после
замены получается система. Напомним, что через $(y_i')_{i=0}^n$ мы обозначаем
последовательность, порождённую~$\SC'$. Заметим, что $y_i'=y_i$
при $i<r$. Положим $t=\ell-r-1\geq 1$.

Обозначим $a=y_{r-2}=y'_{r-2}$, $b=y_{r-1}=y'_{r-1}$; заметим, что $b\leq 2a$
по предложению~\ref{vdvoe}. Тогда $y_r=2b$, и по предложению~\ref{ryad} получаем
$$
  y_{\ell-1}=y_{r+t}=F(t)y_r+F(t-1)y_{r-1}=(2F(t)+F(t-1))b=F(t+2)b.
$$
Аналогично имеем
\begin{align*}
  y'_{\ell-2}=y'_{r+t-1}&=F(t)y_{r-1}+F(t-1)y_{r-2}=F(t)b+F(t-1)a\geq \\
    &\geq F(t)b+\frac{F(t-1)}2b=\frac{F(t)}2b,\\
  y'_{\ell-1}=2y'_{\ell-2}&\geq F(t+2)b=y_{\ell-1}.
\end{align*}

Далее, $y_\ell=y_{\ell-1}+y_{\theta(\ell)}$,
$y_{\ell}'=y_{\ell-1}'+y_{\theta'(\ell)}'$. При этом, если $\ell=\pi(r)$, то
$y_{\theta'(\ell)}'=y_{\ell-3}'\geq y_{r-2}'=y_{r-2}=y_{\theta(r)}$. Иначе
$\theta'(\ell)=\theta(\ell)$, и либо $\theta(\ell)=\ell-1$, либо
$\theta(\ell)\leq r-1$, ибо индексы $r+1,\dots,\ell-1$ рядовые. В любом случае
получаем $y_{\theta'(\ell)}'\geq y_{\theta(\ell)}$, а потому и $y_\ell'\geq
y_\ell$.

Мы готовы проверить, что условия леммы~\ref{uluchsh} об улучшении выполнены, откуда будет следовать требуемое. Условие~\itref{2} уже проверено; условие~\itref{1} очевидно. Условие же~\itref{3} очевидно для всех $i\notin\{\psi(r),d^{-1}(r),d^{-1}(\ell)\}$, ибо тогда из $\theta(i)<\ell-1\leq i-2$ следует $\theta(i)=\theta'(i)\leq r-1$ и $y_{\theta(i)}=y'_{\theta'(i)}$. Рассмотрим оставшиеся случаи. При $i=\psi(r)=\psi'(\ell-1)$ имеем $y_{\theta'(i)}'=y_{\ell-2}'\geq y_{r-1}'=y_{r-1}=y_{\theta(i)}$. При $i=d^{-1}(r)$ имеем $y_{\theta'(i)}'=y_{\ell-3}'\geq y_{r-2}'=y_{r-2}=y_{\theta(i)}$. Наконец, случай $i=d^{-1}(\ell)$ невозможен, ибо $\ell$ --- неплохой, т.е. $\ell\notin d(K)$.
\qed

\note Подобную операция {\em замены} одного индекса другим мы будем
описывать многократно. В дальнейшем мы не будем описывать соответствующую
замену функций $\psi$, $\pi$, считая её подразумевающейся.

\begin{prop}[перемена мест]
  \label{ch-pl}
  Пусть $2\leq r\leq n-2$, $r\in I$, $r+1\in J\cup K$, причём $r+1\notin\{\psi(r),\pi(r)\}$,
  а индекс $r+2$ --- неплохой. Если $r+2$ --- регулярный, то заменим в~$I$
  индекс $r$ на $r+2$, а в~$J$ или в~$K$ --- индекс~$r+1$ на $r$. Иначе заменим в $I$
  индекс $r$ на $r+1$, а в~$J$ или в~$K$ --- индекс $r+1$ на $r$.

  Тогда полученная система~$\SC'$ улучшает~$\SC$.
\end{prop}

\proof
Заметим сразу, что $y_i=y_i'$ при $i<r$. Кроме того, условие
$r+1\notin\{\psi(r),\pi(r)\}$ гарантирует, что после замены получается система.
Обозначим $b=y_{r-1}$, $p=y_{\theta(r+1)}$.
Возможны три случая.

1. Пусть индекс~$r+2$ --- регулярный. Тогда
\begin{align*}
  y_r&=2b, & y_{r+1}&=2b+p, & y_{r+2}&=4b+p, \\
  y'_r&=b+p, & y'_{r+1}&=2b+p=y_{r+1}, & y'_{r+2}&=4b+2p\geq y_{r+2}.
\end{align*}
Проверим условия леммы об улучшении при $\ell=r+2$. Условие~\itref{1} очевидно,
а \itref{2} уже проверено. Условие~\itref{3} требует проверки лишь при
$i\in\{\psi(r),d^{-1}(r)\}$ (иначе $\theta(i)=\theta(i')<r$ и
$y_{\theta(i)}=y'_{\theta'(i)}$). Если $i=\psi(r)=\psi'(r+2)$, то
$y'_{\theta'(i)}=y'_{r+1}\geq y'_{r-1}=y_{\theta(i)}$. Если же $r=d(i)$, то
$y'_{\theta'(i)}=y'_r\geq y'_{r-2}=y_{\theta(r)}$.

2. Пусть теперь индекс $r+2$ --- исключительный и неплохой. Тогда
\begin{align*}
  y_r&=2b, & y_{r+1}&=2b+p, & y_{r+2}&=4b+2p, \\
  y'_r&=b+p, & y'_{r+1}&=2b+2p\geq y_{r+1}, & y'_{r+2}&=4b+4p\geq y_{r+2}.
\end{align*}
Опять проверим условия леммы об улучшении при $\ell=r+2$. Условия~\itref{1} и~\itref{2} верны. Условие~\itref{3} требует проверки лишь при $i\in\{\psi(r),d^{-1}(r)\}$ (напомним, что $r+2$ --- неплохой); эта проверка производится аналогично предыдущему случаю.

3. Наконец, пусть индекс~$r+2$ --- штрафной или пеневой. Обозначим $q=y_{\theta(r+2)}$. Тогда
\begin{align*}
  y_r&=2b, & y_{r+1}&=2b+p, & y_{r+2}&=2b+p+q, \\
  y'_r&=b+p, & y'_{r+1}&=2b+2p\geq y_{r+1}, & y'_{r+2}&=2b+2p+q\geq y_{r+2}.
\end{align*}
Проверка условий леммы об улучшении проводится аналогично первому случаю.
\eop

\begin{cor}[о разделении]
  \label{razd}
  Пусть $2\leq p\leq q<n$, причём $q+1$ --- неплохой нерядовой индекс.   Обозначим $T=[p,q]$. Предположим, что $T\cap K=\emp$, и все индексы из множества $I\cap T$, кроме, возможно, наименьшего из них --- неплохие.

  Тогда существует система $\SC'=(n,I',J',K',\psi',\pi')$, улучшающая~$\SC$; при этом $\SC$ и $\SC'$ отличаются лишь на отрезке~$T$ (формально говоря, $I\setminus T=I'\setminus T$, $J\setminus T=J'\setminus T$, $K\setminus T=K'\setminus T$, и $\psi(i)=j\iff \psi'(i)=j$, $\pi(i)=k\iff \pi'(i)=k$ для любых $i,j,k\notin T$), и множество~$I'\cap T$ находится правее, чем~$J'\cap T$. Более того, $|I'\cap T|=|I\cap T|$, $|J'\cap T|=|J\cap T|$, и $K'\cap T=\emp$.
\end{cor}

\proof
Если $I_1=I\cap T$ уже находится правее~$J_1=J\cap T$ (в частности, если одно
из этих множеств пусто), то можно положить $\SC'=\SC$. Иначе выберем
$$
  i_0=\min I_1, \qquad j=\min\{j\in J_1: j>i_0\}, \qquad i=\max\{i\in I_1: i<j\}.
$$
Тогда $i$ можно заменить на~$j-1$ по предложению~\ref{iskl-right}, а затем поменять их местами по одному из вариантов предложения~\ref{ch-pl}. Последнее возможно, так как $j+1$ не может быть плохим индексом по условию, а также $\pi(i)>j$ (ибо $I\cap K=\emp$). Поскольку сумма исключительных индексов строго возрастает, серией таких замен мы рано или поздно добьёмся требуемого.

Осталось заметить, что при каждой замене мощности множеств $I\cap T$, $J\cap T$ и $K\cap T$ не менялись.
\eop

%
Ещё одно преобразование связано только с изменением функции~$\psi$, то есть с
<<переназначением>> штрафных индексов.

\begin{prop}[о переназначении двух штрафов]
  \label{2straf}
  Пусть $i_1<i_2$ --- некоторые исключительные индексы, а $j_1<j_2$ ---
  соответствующие им штрафные (т.е. $\psi(i_s)=j_s$ при $s=1,2$), причём $j_1>\pi(i_2)$. Изменим
  функцию~$\psi$ на элементах $i_1$, $i_2$, полагая $\psi'(i_s)=j_{3-s}$ при
  $s=1,2$. Тогда получилась система $\SC'$, являющаяся улучшением системы~$\SC$.
\end{prop}

\proof
Условия $i_1<i_2<\pi(i_2)<j_1<j_2$ гарантируют, что $\SC'$ --- система.
Заметим, что $y_t=y_t'$ при $t<j_1$. Обозначим $a_s=y_{i_s-1}=y_{\theta(j_s)}=y'_{\theta'(j_{3-s})}$ при
$s=1,2$, $\delta=a_2-a_1\geq 0$. Тогда $y_{j_1}'=y_{j_1-1}+a_2=y_{j_1}+\delta$.
Непосредственная индукция показывает, что $y_t'\geq y_t+\delta$ при $j_1\leq
t<j_2$. Тогда $y_{j_2}'=y_{j_2-1}'+a_1\geq y_{j_2-1}+a_1+\delta=y_{j_2}$. Тогда
нетрудно видеть, что все условия леммы~\ref{uluchsh} об улучшении при $\ell=j_2$ выполнены.
\eop

\begin{cor}[о переназначении штрафов]
  \label{mstraf}
  Пусть $k\geq 2$, $i_1<\dots<i_k$ --- некоторые исключительные индексы, а $j_1>\dots>j_k$ ---
  штрафные индексы, причём $\psi(\{i_1,\dots,i_k\})=\{j_1,\dots,j_k\}$. Изменим
  функцию~$\psi$ на элементах $i_1,\dots,i_k$, полагая $\psi'(i_s)=j_s$ при
  $s=1,\dots,k$. Тогда, если $\SC'=(n,I,J,K,\psi',\pi)$ --- система, то $\SC'$ ---
  улучшение системы~$\SC$.
\end{cor}

\proof
Индукция по~$k$. При $k=2$ это --- предыдущее предложение. Пусть $k>2$. Если $\psi(i_k)=j_k$, то можно непосредственно применить предположение индукции. Пусть теперь $\psi(i_k)=j_s$ при $s<k$; тогда $j_k=\psi(i_t)$ при некотором~$t<k$. Применим предложение~\ref{2straf} к индексам $i_t$, $i_k$, $j_k$, $j_s$. Поскольку $\SC'$ --- система, то $i_\ell<i_k<\pi(i_k)<j_k<j_s$, поэтому при замене функции $\psi$ переназначением $\psi''(i_t)=j_s>j_k$, $\psi''(i_k)=j_k$ также получается система~$\SC''$, являющаяся улучшением~$\SC$. Для неё опять можно применить предположение индукции, ибо $\psi''(i_k)=j_k$.
\eop

\note Для того, чтобы в условиях следствия~\ref{mstraf} $\SC'$ оказалась системой, достаточно, например, чтобы
выполнялось условие $|\pi(\{i_1,\dots,i_k\})|=1$.

\subsection{Случай единственной пени}

Разберём сначала случай, когда $|K|=1$. В этом случае оказывается верна следующая лемма.

\begin{lm}
  Пусть последовательности $(y_i)$ и $(x_i)$ порождены системами $\SC=(n,I,J,K,\psi,\pi)$ и~$\OC_n$ оответственно, причём  $2\in I$, $|K|=1$ и $x_n\geq y_n$. Тогда для любых допустимых последовательностей $\Pi_{\OC_n}(a,b)=(x_i')$ и $\Pi_{\SC}(a,b)=(y_i')$ имеем $x_n'\geq y_n'$.
  \label{compare-emp}
\end{lm}

\proof
Обозначим $(a_i)=\Pi_{\OC_n}(1,0)$, $(b_i)=\Pi_\SC(1,0)$; пусть $K=\{k\}$. Заметим, что $b_i=0$ при $i<k$, поскольку при всех таких индексах $\theta(i)>0$. Кроме того, поскольку $k\geq 3$, мы имеем $a_k\geq 1=b_k$ и $a_{k+1}\geq 2=b_k+b_0=b_k+b_{\theta(k)}=b_{k+1}$.

Далее, никакой индекс $i\geq k$ --- не исключительный. Покажем индукцией по~$i\geq k$, что $b_i\leq a_i$. Действительно, при $i=k$ и $i=k+1$ утверждение уже доказано; если же $i\geq k+2$, то $b_i=b_{i-1}+b_{\theta(i-1)}$; если $\theta(i-1)\neq i-2$, то $b_{\theta(i-1)}=0$ и $b_i= b_{i-1}\leq a_{i-1}\leq a_i$; иначе $b_i\leq b_{i-1}+b_{i-2}=a_{i-1}+a_{i-2}=a_i$, что и требовалось.

Итак, мы получаем, что $b_n\leq a_n$.
Наконец, заметим, что $\Pi_{\OC_n}(a,b)=(b/2)\Pi_{\OC_n}(1,2)+(a-b/2)\Pi_{\OC_n}(1,0)$ и $\Pi_\SC(a,b)=(b/2)\Pi_\SC(1,2)+(a-b/2)\Pi_\SC(1,0)$, причём $a-b/2\geq 0$, поскольку эти последовательности допустимы. Значит,
$$
  x_n'=(b/2)x_n+(a-b/2)a_n\geq (b/2)y_n+(a-b/2)b_n=y_n',
$$
что и требовалось доказать.
\qed

\begin{df}
  Для числового множества $X$ определим его {\em сдвиг влево} как $X^-=\{x-1:x\in X\}$. Для числовой функции $\phi$ определим её {\em сдвиг влево} формулой $\phi^-(x)=\phi(x+1)-1$.

  Пусть $\SC=(n,I,J,K,\psi,\pi)$ --- система, в которой $2$ --- регулярный индекс. Определим её {\em сдвиг влево} как систему $\SC^-=(n-1,I^-,J^-,K^-,\psi^-,\pi^-)$.
\end{df}

\begin{lm}
  Пусть $\SC=(n,I,J,K,\psi,\pi)$ --- система с $|K|=1$. Пусть $(x_i)=\Pi_\SC(a,b)$ и $(y_i)=\Pi_{\OC_n}(a,b)$ --- допустимые последовательности. Тогда $x_n\geq y_n$.
  \label{gen-pen1}
\end{lm}

\proof
Предположим противное; выберем из всех допустимых по\-с\-ле\-до\-ва\-тель\-но\-с\-тей~$\Pi_\SC(a,b)$ (при всевозможных $\SC$, $a$ и $b$), противоречащих лемме, ту, для которой $n$ минимально, а из таких --- ту, для которой минимально~$|I|$. Согласно лемме~\ref{compare-emp}, можно считать, что $a=1$, $b=2$. Если $2\notin I$, то $2$ --- регулярный индекс. Значит, для последовательностей $(x_{i+1})_{i=0}^{n-1}=\Pi_{\SC^-}(b,a+b)$ и $(y_{i+1})_{i=0}^{n-1}=\Pi_{\OC_{n-1}}(b,a+b)$ утверждение леммы верно, то есть $x_n\geq y_n$, что не так. Итак, $2\in I$.

Пусть $K=\{k\}$. Предположим, что $I\neq \{2,3,\dots,k-1\}$. Положим $i_0=\min\{i\geq 2: i\notin I\}$, $j_0=\min\{\ell\in L:\ell>i_0\}$, $t=j_0-i_0$. По лемме~\ref{iskl-right}, можно последовательно сдвинуть все индексы~$i_0-1, i_0-2,\dots,2$ в индексы $i_0+t-1,\dots,2+t$ соответственно, улучшив систему~$\SC$. При этом $2$ не является исключительным для новой системы, что невозможно. Значит, $I=\{2,\dots,k-1\}$. Далее, согласно следствию~\ref{mstraf} о переназначении штрафов, можно считать, что $\psi(2)>\psi(3)>\dots>\psi(k-1)$.

Предположим, что $|I|\geq 2$. Положим $I'=(I\setminus\{2\})^-$, $J'=(J\setminus\{\psi(2)\})^-$, $K'=\{k-1\}$, $\psi'=\psi^-\big|_{I'}$, $\pi'(i)=k-1$ для всех $i=2,\dots,k-2$. Тогда $\SC'=(n,I',J',K',\psi',\pi')$ --- система; пусть $(y_i')=\Pi_{\SC'}$. Покажем, что $y_n'\geq y_n$; это будет противоречить исходному выбору, ибо $|I'|< |I|$.

Положим $t=\max J$, $t'=\max J'$; тогда $d=t-t'\geq 2$. Положим $\SC''=(t',I',J',K',\psi',\pi')$. Пусть $(a_i)=\Pi_{\SC''}=(y_i')_{i=0}^{t'}$. Положим $p=a_{t'-1}$, $q=a_{t'}$. Заметим, что $p\geq 5$, $q=a_{t'-1}+a_1=p+2\geq 7$. Тогда по предложению~\ref{ryad},
$y_{t-1}'=F(d-1)q+F(d-2)p$, $y_t'=F(d)q+F(d-1)p$. С другой стороны, имеем $y_1=2$, $y_2=4$; значит, отрезок последовательности $(y_{i+1})_{i=0}^{t'}$ строится так же, как и $\Pi_{\SC''}(2,4)=(2a_i)$, за единственным исключением: $y_k=y_{k-1}+1$, в то время как $2a_{k-1}=2a_{k-2}+2$. Тогда нетрудно видеть, что $y_{t'}\leq 2p$, $y_{t'+1}\leq 2q$. Теперь, снова по предложению~\ref{ryad}, получаем
$y_{t-1}\leq 2F(d-2)q+2F(d-3)p$. Тогда
$$
  y_t=y_{t-1}+y_1\leq 2F(d-2)q+2F(d-3)p+2\leq F(d)q+F(d-1)p= y_t',
$$
поскольку $F(d)\geq 2F(d-2)$, $F(d-1)\geq 2F(d-3)$, причём хотя бы одно из этих неравенств строгое. Далее, пусть $t<n$; покажем тогда, что $y_{t+1}'\geq y_{t+1}$. Поскольку $q=p+2$, имеем
\begin{align*}
  y_{t+1}'-y_{t+1}&=(y_t'+y_{t-1}')-(y_t+y_{t-1})
    \geq (F(d+1)q+F(d)p)-(4F(d-2)q+4F(d-3)p+2)=\\
  &=(F(d+2)-4F(d-1))p+2(F(d+1)-4F(d-2)-1)=\\
  &=(2F(d-2)-F(d-1))(p-2)+2(2F(d-1)-F(d)-1).
\end{align*}
Поскольку $2F(d-2)\geq F(d-1)$ и $2F(d-1)\geq F(d)$, причём хотя бы одно из этих неравенств строгое, получаем $y_{t+1}'-y_{t+1}\geq \min\{p-4,0\}=0$. Итак, $y_{t+1}\leq y'_{t+1}$, $y_t\leq y'_t$, откуда и следует, что $y_n\leq y_n'$.

Наконец, пусть $|I|=1$. Тогда аналогично $\SC$ заменяется на пустую систему с увеличением последнего члена.
\qed

\subsection{Общая оценка}

Теперь мы готовы к доказательству общей оценки. Сначала докажем лемму,
позволяющую сделать ключевой индукционный шаг с применением
леммы~\ref{gen-pen1}.

Для каждого~$k\in K$ определим его {\em отрезок влияния} $A(k)=[d(k),\max
\psi(\pi^{-1}(k))]$. Иными словами, отрезок влияния пеневого индекса~$k$ --- это
минимальный отрезок, содержащий все исключительные и штрафные индексы,
соответствующие~$k$. Назовём индекс $k\in K$ {\em выделенным}, если на отрезке
$A(k)$ нет индексов, соответствующих другому пеневому индексу (иначе говоря,
$\pi(I\cap A(k))=\pi(\psi^{-1}(J\cap A(k)))=K\cap A(k)=\{k\}$).

\begin{lm}[о выделении
]
  \label{vydel}
  Пусть $|K|\geq 2$. Тогда существует система~$\SC'$, улучшающая~$\SC$, такая, что
  в ней $|K'|\leq |K|$, причём либо в~$\SC'$ существует выделенный пеневой индекс~$k_0$, либо $|K'|<|K|$. При этом в первом случае имеем $d(K)\cap[k_0,n]=\emp$.
\end{lm}

\proof
Процесс улучшения будет проходить в несколько шагов. На каждом шаге мы будем
улучшать систему, добиваясь выполнения некоторого свойства (своего для каждого шага). При этом это свойство
будет сохраняться и при всех последующих шагах.

Пусть $i_0=\max d(K)$,  $k_0=\pi(i_0)$.

{\em Шаг 1.} Предположим, что существует пеневой индекс~$k$, лежащий на интервале
$(i_0,k_0)$. Переопределим функцию~$\pi$, полагая $\pi'(i)=k$ для всех
$i\in[i_0,k]\cap \pi^{-1}(k_0)$. Из определения~$i_0$ следует, что значение
$d(k)$ не изменилось. Если после этого оказалось, что $\pi'^{-1}(k_0)=\emp$, то
выкинем $k_0$ из~$K$. Нетрудно видеть, что получилась система, причём являющаяся улучшением исходной (все члены порождённой последовательности вплоть до $k_0$-го не изменились,
последующие не уменьшились); при этом, если $|K'|=|K|$, то значение $d(k_0)$ увеличилось. После нескольких таких шагов мы либо уменьшим~$|K|$ (тем самым добившись требуемого), либо добьёмся того, что $(i_0,k_0)\cap K=\emp$.

Итого, можно считать, что $(i_0,k_0)\cap K=\emp$.

{\em Шаг 2.} Переопределим функцию $\pi$, полагая $\pi'(i)=k_0$ для всех $i\in
[i_0,k_0]\cap I$. Получилась снова система, причём, поскольку множество
$(i_0,k_0]\cap I$ не содержит плохих индексов, порождённая последовательность не
изменилась.

Итого, можно считать, что $\pi([i_0,k_0]\cap I)=k_0$.

{\em Шаг 3.} Рассмотрим отрезок~$T=[i_0,k_0-1]$; на нём нет пеневых индексов, а
плохим является только индекс $i_0$. По следствию~\ref{razd} о разделении можно
так перестроить систему~$\SC$ на отрезке~$T$, что на отрезке~$T$ все штрафные
индексы будут левее всех исключительных. При этом значение $i_0$ может только увеличиться, а свойство шага~2 сохраняется.

Иными словами, можно считать, что $[i_0,k_0]\cap J=\emp$.

{\em Шаг 4.} Обозначим $I_1=[i_0,k_0]\cap I=\{i_0,\dots,i_t\}$ ($i_0<\dots<i_t$), $J_1=\psi(I_1)$. По
следствию~\ref{mstraf} о переназначении штрафов, можно считать, что
$\psi(i_s)>\psi(i_r)$ при $0\leq s<r\leq t$. Обозначим $j_s=\psi(i_s)$ при
$0\leq s\leq t$.

{\em Шаг 5.} Предположим, что на отрезке $[i_0,j_0]$ содержится ещё какой-то
пеневой индекс~$k\neq k_0$ (тогда $k>k_0$; мы выбираем $k$ наименьшим
возможным). Пусть $s_0=\max\{s: j_s>k\}$. Тогда можно переопределить
функцию~$\pi$ на элементах $i_0,\dots,i_{s_0}$, полагая $\pi'(i_s)=k$ при
$0\leq s\leq s_0$ (при этом, если $s_0=t$, то надо ещё выкинуть $k_0$ из~$K'$,
уменьшив тем самым~$|K|$; в противном случае будем иметь $d'(k_0)=i_{s_0+1}$).
При этом значение $d(k)$ не изменится по выбору~$i_0$. Нетрудно видеть, что
получилась система, улучшающая исходную: члены порождённой последовательности
вплоть до~$k_0$-го не изменились, а дальнейшие не уменьшились.

При этом в изменённой системе (для новых значений $i_0$, $k_0$) выполнено
условие $[i_0,j_0]\cap K=\{k_0\}$.

{\em Шаг 6.} Пусть теперь $J_1=J\cap [k_0,j_0]$.
Покажем, что можно переназначить штрафы так, чтобы индексы~$j_0,\dots,j_t$
являлись минимальными индексами в $J_1$ (иначе говоря, чтобы
$\{j_0,\dots,j_t\}=J\cap [k_0,j_0]$). Пусть это не так, то есть для некоторого
$0\leq s\leq t$ существует $j\in J\setminus\{j_0,\dots,j_t\}$ такой, что
$j<j_s$; пусть $i=\psi^{-1}(j)<i_0$. Можно считать, что $s$ --- максимальный
индекс с этим свойством, а~$j$ --- минимальный для этого~$s$. Тогда по
предложению~\ref{2straf} о переназначении двух штрафов можно переназначить
$\psi'(i)=j_s$, $\psi'(i_s)=j$; ясно, что получится система, причём в ней для
индекса~$s$ уже не будет существовать таких~$j$. Повторяя процедуру, в конце
концов добьёмся требуемого. Заметим, что в процессе переназначений порядок
элементов $\psi(i_0),\dots,\psi(i_t)$ остаётся неизменным.

Итак, мы добились того, что $\{j_0,\dots,j_t\}=[i_0,j_0]\cap J$ (напомним, что
$[i_0,k_0]\cap J=\emp$ уже после Шага~3).

{\em Шаг 7.} Предположим, наконец, что $[i_0,j_0]\cap L\neq \{i_0,\dots,i_t,j_0,\dots,j_t,k_0\}$;
по результатам предыдущих шагов, <<лишними>> элементами могут быть только
исключительные индексы, лежащие на отрезке $[k_0,j_0]$. Тогда пеневые индексы,
им соответствующие, больше~$j_0$. Положим $f_0=\min\{f\in L: f>j_0\}$. Из
выбора~$i_0$ следует, что $f_0$ --- неплохой. Тогда по следствию~\ref{razd} о
разделении, существует улучшение~$\SC'$ нашей системы, отличающееся от неё лишь
на отрезке~$[k_0+1,f_0-1]$, в котором уже $\ds\left[k_0,\max_{0\leq s\leq t}
\psi(i_s)\right]\cap I=\emp$.

Суммируя предыдущие результаты, видим, что в полученной системе выполняется
соотношение
$$
  [i_0,j_0]\cap L=\{i_0,\dots,i_t,j_0,\dots,j_t,k_0\}.
$$
Таким образом, индекс~$k_0$ в ней является выделенным. Кроме того, по определению~$i_0$, мы имеем $d(K)\cap [k_0,n]=\emp$.
\eop

\begin{thm}
  \label{y<F}
  $y_n\leq F(n+1)$.
\end{thm}

\proof
Индукция по $|K|$. При $|K|=0$ доказывать нечего, при $|K|=1$ утверждение
следует из леммы~\ref{gen-pen1}.

Пусть $|K|\geq 2$. Применяя лемму~\ref{vydel} о выделении, мы либо уменьшим
$|K|$ (после чего применимо предположение индукции), либо получим систему с тем
же значением $|K|$, в которой некоторый пеневой индекс~$k_0$ --- выделенный.
Полученную систему опять будем обозначать через~$\SC$.

Итого, пусть индекс~$k_0\in K$ выделен. Пусть $\SC'=(n,I',J',K',\psi',\pi')$, где
$$
  K'=K\setminus\{k_0\}, \qquad I'=I\setminus \pi^{-1}(k_0), \qquad
  J'=J\setminus \psi(\pi^{-1}(k_0)),
  \qquad \psi'=\psi|_{I'}, \qquad \pi'=\pi|_{I'};
$$
нетрудно видеть, что $\SC'$ --- система. Пусть $(y_i')$ --- последовательность, порождённая~$\SC'$. Поскольку~$|K'|=|K|-1$, по предположению индукции
$y_n'\leq F(n+1)$. Для завершения доказательства достаточно доказать, что $\SC'$
улучшает~$\SC$.

Положим
\begin{gather*}
  I_0=\pi^{-1}(k_0), \quad i_0=\min I_0, \qquad
  J_0=\psi(I_0), \quad j_0=\max J_0.
\end{gather*}
Заметим, что $y_i'=y_i$ при $i<i_0$.
Пусть $j$ --- максимальный индекс такой, что все индексы из полуинтервала $(j_0,j]$
--- рядовые (таким образом, если $j_0=\max L$, то $j=n$, иначе $j=\min\{j\in L: j>j_0\}-1$).

Неформально говоря, поскольку индекс~$k_0$ "--- выделенный,
последовательность $(y_i)$ ведёт себя на отрезке $[i_0,j]$ в точности как последовательность, допустимая для <<сдвинутой>> системы $(j-i_0+2,I_0,J_0,K_0,\psi|_{I_0},\pi|_{I_0})$. Формализуем это утверждение.

Положим $i_*=i_0-2$. Определим систему~$\SC_1=(j-i_*,I_1,J_1,K_1,\psi_1,\pi_1)$
следующим образом:
\begin{gather*}
  I_1=I_0-i_*, \qquad J_1=J_0-i_*, \qquad K_1=\{k_0-i_*\}, \\
  \psi_1(i-i_*)=\psi(i)-i_*, \qquad \pi_1(i-i_*)=\pi(i)-i_*.
\end{gather*}
Пусть $(x_i)_{i=0}^{j-i_*}=\Pi_{\SC_1}(y_{i_*},y_{i_*+1})$. Тогда
непосредственная индукция показывает, что $x_i=y_{i+i_*}$ при всех $0\leq i\leq
j-i_*$, ибо оба члена получаются из предыдущих по одинаковым правилам.
Аналогично, если $(x_i')_{i=0}^{j-i_*}$ --- допустимая последовательность для пустой системы с теми же начальными условиями, то $x_i'=y_{i+i_*}'$ при всех $0\leq i\leq j-i_*$. По следствию~\ref{gen-pen1}, имеем теперь $y_j=x_{j-i_*}\leq x_{j-i_*}'=y_j'$. Если $j=n$, то это и есть
требуемое неравенство.

Пусть, наконец, $j<n$. Покажем, что последовательности $(y_i)$ и $(y_i')$
удовлетворяют условиям леммы~\ref{uluchsh} об улучшении при $\ell=j+1$. Заметим, что
$y_{j+1}=y_j+y_{\theta(j+1)}$, $y_{j+1}'=y_j'+y_{\theta'(j+1)}'$, причём индекс
$\theta=\theta(j+1)=\theta'(j+1)$ либо равен $j$ (если $j+1\in I$), либо не
превосходит $i_0-1$. Значит, $y_{\theta(j+1)}\leq y'_{\theta(j+1)}$, откуда
следуют условия~\itref{2}. Напомним, что из утверждения леммы~\ref{vydel}
вытекает, что индекс~$j+1$ неплохой.

Далее, при любом $i\geq j+1$ имеем $\theta=\theta(i)=\theta(i)'$, что
доказывает~\itref{1}. Кроме того, при этих же значениях~$i$ либо $\theta(i)\geq
j$, либо $\theta(i)\leq i_0-1$ (это следует из того, что $j+1$ --- неплохой, а
$k_0$ --- выделенный). Поэтому условие~\itref{3} также выполнено, ибо
$y_s=y_s'$ при $s\leq i_0-1$.

Итак, по лемме об улучшении $y_n\leq y_n'$, что и требовалось доказать.
\eop

\medskip
Теперь мы можем доказать основную теорему.

\textbf{Доказательство теоремы~\ref{main'}.} Пусть $v_1,\dots,v_n$ --- все
развилки в слове~$W$, упорядоченные по неубыванию значимости, $z_i=r(v_i)$. По
предложению~\ref{z<y}, $z_n\leq y_n$ для некоторой порождённой последовательности $(y_i)_{i=0}^n$. По теореме~\ref{y<F}, $y_n\leq F(n+1)$.
Значит, и $|u|=z_n\leq F(n+1)$, что и требовалось доказать.
\eop

\section{Алфавит из произвольного количества букв}


\begin{thebibliography}{99}
\bibitem{Ufn}
В.А. Уфнаровский. Комбинаторные и асимптотические методы в алгебре
// ВИНИТИ, 1990 Сер. Совр. пробл. Математики. Фундаментальные направления. Т.57.
 М.  C.5 --- 177

\bibitem{Kur}
А.Г. Курош.  Проблеммы теории колец, связанные с проблеммой Бернсайда о
периодических группах //Изв. АН СССР сер. мат. 1941. Т.5. C. 233
--- 240

\bibitem{BBL}
А.Я. Белов, В.В. Борисенко, В.Н. Латышев. Мономиальные алгебры
//ВИНИТИ, 2002. Итоги науки и техники. Сер. Современная математика и
ее приложения. Тематические обзоры. Т.26. М. C.35 --- 214



\bibitem{cheln} Г.Р. Челноков. О числе запретов, задающих периодическую
последовательность. // Модел. и анализ информ. систем. Т.~13, \No~3 (2007),
66--70.

\bibitem{lavrov} P. Lavrov. Number of restrictions required for periodic words in finite alphabet. 	 {\tt arXiv:1209.0220.}

\end{thebibliography}
\end{document}